# SPLITTING OF LIFTINGS IN PRODUCTS OF PROBABILITY SPACES[1]

By W. Strauss, N. D. Macheras and K. Musiał

*Universität Stuttgart, University of Piraeus and Wrocław University*

We prove that if $(X, \mathfrak{A}, P)$ is an arbitrary probability space with countably generated $\sigma$-algebra $\mathfrak{A}$, $(Y, \mathfrak{B}, Q)$ is an arbitrary complete probability space with a lifting $\rho$ and $\widehat{R}$ is a complete probability measure on $\mathfrak{A} \widehat{\otimes}_R \mathfrak{B}$ determined by a regular conditional probability $\{S_y : y \in Y\}$ on $\mathfrak{A}$ with respect to $\mathfrak{B}$, then there exist a lifting $\pi$ on $(X \times Y, \mathfrak{A} \widehat{\otimes}_R \mathfrak{B}, \widehat{R})$ and liftings $\sigma_y$ on $(X, \widehat{\mathfrak{A}}_y, \widehat{S}_y)$, $y \in Y$, such that, for every $E \in \mathfrak{A} \widehat{\otimes}_R \mathfrak{B}$ and every $y \in Y$, $[\pi(E)]^y = \sigma_y([\pi(E)]^y)$.
Assuming the absolute continuity of $R$ with respect to $P \otimes Q$, we prove the existence of a regular conditional probability $\{T_y : y \in Y\}$ and liftings $\varpi$ on $(X \times Y, \mathfrak{A} \widehat{\otimes}_R \mathfrak{B}, \widehat{R})$, $\rho'$ on $(Y, \mathfrak{B}, \widehat{Q})$ and $\sigma_y$ on $(X, \widehat{\mathfrak{A}}_y, \widehat{S}_y)$, $y \in Y$, such that, for every $E \in \mathfrak{A} \widehat{\otimes}_R \mathfrak{B}$ and every $y \in Y$,
$[\varpi(E)]^y = \sigma_y([\varpi(E)]^y)$
and
$$\varpi(A \times B) = \bigcup_{y \in \rho'(B)} \sigma_y(A) \times \{y\} \quad \text{if } A \times B \in \mathfrak{A} \times \mathfrak{B}.$$
Both results are generalizations of Musiał, Strauss and Macheras [*Fund. Math.* **166** (2000) 281–303] to the case of measures which are not necessarily products of marginal measures. We prove also that liftings obtained in this paper always convert $\widehat{R}$-measurable stochastic processes into their $R$-measurable modifications.

**1. Preliminaries.** If $(Z, \mathfrak{Z}, S)$ is a probability space, then we denote by $\widehat{\mathfrak{Z}}$ the completion of $\mathfrak{Z}$ with respect to $S$ and by $\widehat{S}$ the completion of $S$. We write $\mathcal{L}^\infty(S) := \mathcal{L}^\infty(Z, \mathfrak{Z}, S)$ for the space of bounded $\mathfrak{Z}$-measurable real-valued functions. Functions equal a.e. are not identified.

We use the notion of *lower density* and *lifting* in the sense of [7]. It is known (cf. [7]) that there is a 1–1 correspondence among liftings on $\mathfrak{Z}$ with

Received July 2002; revised June 2003.

[1]Supported in part by KBN Grant 5 P03A 016 21 and by NATO Grant PST.CLG.977272.

respect to $S$ and liftings on $\mathcal{L}^\infty(S)$. Saying that $\tau$ is a lifting on $(Z, \mathfrak{Z}, S)$, we mean that $\tau$ is a lifting on $\mathfrak{Z}$ and on $\mathcal{L}^\infty(S)$. $\Lambda(S)$ denotes the system of all liftings on $(Z, \mathfrak{Z}, S)$. Similarly, $\vartheta(S)$ is the collection of all lower densities on $(Z, \mathfrak{Z}, S)$.

Throughout what follows, $(X, \mathfrak{A}, P)$ and $(Y, \mathfrak{B}, Q)$ are probability spaces. $\mathfrak{A} \times \mathfrak{B}$ is the product algebra generated by $\mathfrak{A}$ and $\mathfrak{B}$ in $X \times Y$, and $\mathfrak{A} \otimes \mathfrak{B} := \sigma(\mathfrak{A} \times \mathfrak{B})$ is the product $\sigma$-algebra generated by $\mathfrak{A} \times \mathfrak{B}$. $P \otimes Q$ is the product measure on $\mathfrak{A} \otimes \mathfrak{B}$, $\mathfrak{A} \widehat{\otimes} \mathfrak{B}$ is the completion of $\mathfrak{A} \otimes \mathfrak{B}$ with respect to $P \otimes Q$ and $P \widehat{\otimes} Q$ is the completion of $P \otimes Q$. We write $\mathfrak{A}_0 = \{A \in \mathfrak{A} : P(A) = 0\}$, $\mathfrak{B}_0 = \{B \in \mathfrak{B} : Q(B) = 0\}$ and $\widehat{\mathfrak{B}}_0 = \{B \subset Y : Q^*(B) = 0\}$, where $Q^*$ is the outer measure generated by $Q$.

$R$ is always a probability measure on $\mathfrak{A} \otimes \mathfrak{B}$, such that $P$ and $Q$ are the marginals of $R$. By $(X \times Y, \mathfrak{A} \widehat{\otimes}_R \mathfrak{B}, \widehat{R})$ we denote the completion of the probability space $(X \times Y, \mathfrak{A} \otimes \mathfrak{B}, R)$. $E_{\mathfrak{C}}(f)$ denotes a version of the conditional expectation of a function $f \in \mathcal{L}^\infty(P)$ with respect to the $\sigma$-algebra $\mathfrak{C} \subset \mathfrak{B}$. An element $C \neq \varnothing$ of an algebra $\mathfrak{C}$ is an atom of $\mathfrak{C}$ if it cannot be decomposed into two disjoint nonempty elements of $\mathfrak{C}$. It follows from the context whether we consider the atoms of an algebra or the atoms of a measure defined on that algebra.

DEFINITION 1.1. Assume that for every $y \in Y$ there is a probability $S_y$ on $\mathfrak{A}$ such that:

(D1) for every $A \in \mathfrak{A}$, the map $y \to S_y(A)$ is $\mathfrak{B}$-measurable;
(D2) $R(A \times B) = \int_B S_y(A)\, dQ(y)$ for all $A \in \mathfrak{A}$ and all $B \in \mathfrak{B}$.

Then the family $\{S_y : y \in Y\}$ is called a *product regular conditional probability* (product r.c.p. for short) on $\mathfrak{A}$ for $R$ with respect to $Q$. One can easily see that the existence of such a product r.c.p. is equivalent to the existence of the classical r.c.p. on the $\sigma$-algebra $\mathfrak{A} \times Y := \sigma(\{A \times Y : A \in \mathfrak{A}\})$ of cylinders based in $\mathfrak{A}$ with respect to the $\sigma$-algebra $X \times \mathfrak{B} := \sigma(\{X \times B : B \in \mathfrak{B}\})$ of cylinders based in $\mathfrak{B}$, on the measure space $(X \times Y, \mathfrak{A} \otimes \mathfrak{B}, R)$. We could use the name of disintegration instead, but it seems that it is better to reserve that term to the general case when $\mathfrak{A}$ is not necessarily countably generated and $S_y$'s may be defined on different domains (cf. [11]). Throughout, we assume that a product regular conditional probability $\{S_y : y \in Y\}$ on $\mathfrak{A}$ with respect to $\mathfrak{B}$ exists. But if $\mathfrak{A}$ is countably generated and $P$ is perfect (cf. [12] for definition), then such an assumption is superfluous since a product r.c.p. always exists (cf. [1], [5] or [11]). A product r.c.p. $\{S_y : y \in Y\}$ on $\mathfrak{A}$ with respect to $\mathfrak{B}$ is said to be *absolutely continuous* with respect to $P$, if $S_y \ll P$ for every $y \in Y$. One can easily see that $R \ll P \otimes Q$ if and only if there exists a product r.c.p. $\{S_y : y \in Y\}$ of $R$ with respect to $Q$ such that $S_y \ll P$ for all $y \in Y$.

The completion of $\mathfrak{A}$ with respect to $S_y$ is denoted by $\widehat{\mathfrak{A}}_y$. The collection of all members of $\mathfrak{A}$ satisfying the relation $S_y(A) = 0$ is denoted by $\mathfrak{A}_{y0}$. If $f \in \mathcal{L}^\infty(S_y) := \mathcal{L}^\infty(X, \mathfrak{A}, S_y)$, then $E_{\mathfrak{C}}^y(f)$ denotes a version of the conditional expectation of $f$ with respect to the $\sigma$-algebra $\mathfrak{C} \subset \mathfrak{A}$.

The whole paper consists of two independent parts. The first one is a continuation of [10]. We have proven in [10] for complete probability spaces $(X, \mathfrak{A}, P)$ and $(Y, \mathfrak{B}, Q)$ that, for given lifting $\rho \in \Lambda(Q)$, there exist liftings $\sigma \in \Lambda(P)$ and $\pi \in \Lambda(P \widehat{\otimes} Q)$ such that

(1) $$\pi(A \times B) = \sigma(A) \times \rho(B) \qquad \text{if } A \times B \in \mathfrak{A} \times \mathfrak{B}$$

and

(2) $$[\pi(E)]^y = \sigma([\pi(E)]^y) \qquad \text{if } E \in \mathfrak{A} \widehat{\otimes} \mathfrak{B} \text{ and } y \in Y.$$

In fact, we have proven more. Namely, the following result holds true: Let $R$ be a measure on $\mathfrak{A} \otimes \mathfrak{B}$ with marginals $P$ and $Q$ and let $\{S_y : y \in Y\}$ be a product r.c.p. of $R$ on $\mathfrak{A}$ with respect to $Q$. If every $S_y$ is equivalent to $P$ (in the sense of absolute continuity), then for given lifting $\rho \in \Lambda(Q)$ there exist liftings $\sigma \in \Lambda(P)$ and $\pi \in \Lambda(\widehat{R})$ such that (1) and (2) hold true. These are in fact the same liftings which were chosen for the product measure $P \otimes Q$, as equivalent measures have the same liftings. (Notice that the properties of the product r.c.p. yield the equivalence of $R$ and $P \otimes Q$.)

Looking at the above result, one may ask immediately what happens if $S_y$'s are not equivalent. Clearly, (1) and (2) may be then senseless, as nonequivalent $S_y$'s have different liftings. So one may ask if for some $\sigma_y \in \Lambda(\widehat{S}_y)$, $y \in Y$, and $\pi \in \Lambda(\widehat{R})$, the equations can have the following form:

(RF) $$\pi(A \times B) = \bigcup_{y \in \rho(B)} \sigma_y(A) \times \{y\} \qquad \text{if } A \times B \in \mathfrak{A} \times \mathfrak{B}$$

and

(SP) $$[\pi(E)]^y = \sigma_y([\pi(E)]^y) \qquad \text{if } E \in \mathfrak{A} \widehat{\otimes}_R \mathfrak{B} \text{ and } y \in Y.$$

One can, however, easily conclude that satisfying the rectangle formula (RF) for all $A \in \mathfrak{A}$ and $B = Y$ yields the absolute continuity of every $S_y$ with respect to $P$. Even more, (RF) forces the following condition (IT) to be satisfied by the product r.c.p. $\{S_y : y \in Y\}$ and by the lifting $\rho \in \Lambda(Q)$:

(IT) $$R(A \times B) = 0 \implies Q(B) = 0 \quad \text{or} \quad S_y(A) = 0 \qquad \text{for all } y \in \rho(B).$$

Clearly, (IT) also yields the absolute continuity of $\{S_y : y \in Y\}$ with respect to $P$. Moreover, we will see in Proposition 2.7 that the conditions (RF) and (IT) are equivalent.

The main result of the first part (Theorem 2.3) shows that if condition (IT) holds true, then (RF) and (SP) are satisfied with properly chosen liftings. If the weaker assumption $R \ll P \otimes Q$ is satisfied, then one can always modify the lifting $\rho$ and the product r.c.p. $\{S_y : y \in Y\}$ to enforce conditions (IT), (RF) and (SP) (see Theorem 2.6).

The second part of the paper deals with quite arbitrary measures $R$, but we assume that the $\sigma$-algebra $\mathfrak{A}$ is separable in the Frechet–Nikodym pseudometric [i.e., there is a countable collection of sets $A_n \in \mathfrak{A}$ such that, given $\varepsilon > 0$ and $A \in \mathfrak{A}$, there is $A_n$ with $P(A \triangle A_n) < \varepsilon$]. We show that in this case one can achieve (SP) (see Theorem 3.6). The proof is completely independent of [10].

**2. Absolutely continuous measures.** Throughout this section, probability spaces $(X, \mathfrak{A}, P)$ and $(Y, \mathfrak{B}, Q)$ are arbitrary but $(Y, \mathfrak{B}, Q)$ is assumed to be complete. According to [10], Theorem 2.9, for given density $\rho \in \vartheta(Q)$, there exist densities $\sigma \in \vartheta(P)$ and $\varphi \in \vartheta(P \widehat{\otimes} Q)$ such that, for all $E \in \mathfrak{A} \widehat{\otimes} \mathfrak{B}$, we have:

(a) $P([\varphi(E)]^y \cup [\varphi(E^c)]^y) = 1$ for all $y \in Y$;
(b) $[\varphi(E)]^y = \sigma([\varphi(E)]^y)$ for all $y \in Y$;
(c) $[\varphi(E)]_x \in \mathfrak{B}$ for all $x \in X$;
(d) $\varphi(A \times B) \supseteq \sigma(A) \times \rho(B)$ whenever $A \in \mathfrak{A}$ and $B \in \mathfrak{B}$.

LEMMA 2.1. *Let* (IT) *be satisfied by a product r.c.p.* $\{S_y : y \in Y\}$ *on* $\mathfrak{A}$ *for $R$ with respect to $Q$ and the density $\rho \in \vartheta(Q)$. If $\tau_y \in \Lambda(\widehat{S}_y)$ for $y \in Y$ is chosen, then there exists a Boolean homomorphism $\varphi_y : \mathfrak{A} \widehat{\otimes}_R \mathfrak{B} \to \tau_y(\widehat{\mathfrak{A}}_y)$ satisfying $\varphi_y(A \times Y) = \tau_y(A)$ for all $y \in Y$ and $\varphi_y(X \times B) = X$ if $y \in \rho(B)$ and $\varphi_y(X \times B) = \varnothing$ otherwise, and such that $\varphi_y(E) = \varnothing$ for all $E \in \mathfrak{A} \widehat{\otimes}_R \mathfrak{B}$ with $\widehat{\widehat{R}}(E) = 0$.*

PROOF. For each $y \in Y$, set $\varphi_y(A \times Y) := \tau_y(A)$, $\varphi_y(X \times B) := X$ if $y \in \rho(B)$ and $\varphi_y(X \times B) := \varnothing$ otherwise. Then set $\varphi_y^0(A \times B) := \varphi_y(A \times Y) \cap \varphi_y(X \times B)$ for $A \in \mathfrak{A}$ and $B \in \mathfrak{B}$. Then $\varphi_y^0 : \mathfrak{A} \times \mathfrak{B} \to \tau_y(\widehat{\mathfrak{A}}_y)$ is a Boolean homomorphism. (IT) simply says that if $R(A \times B) = 0$, then $\varphi_y^0(A \times B) = \varnothing$ for all $y \in Y$, $A \in \mathfrak{A}$ and $B \in \mathfrak{B}$.

Writing $\gamma : \mathfrak{A} \times \mathfrak{B} \to (\mathfrak{A} \times \mathfrak{B})/R$ for the canonical surjection, we find a Boolean homomorphism $\overline{\varphi}_y : (\mathfrak{A} \times \mathfrak{B})/R \to \tau_y(\widehat{\mathfrak{A}}_y)$ such that $\varphi_y^0 = \overline{\varphi}_y \circ \gamma$. Note that $\tau_y(\widehat{\mathfrak{A}}_y)$ is a complete Boolean algebra due to Maharam's theorem (see, e.g., [7], Theorem 3, page 40, or [15], Theorem 3.9, page 1146). For this reason we can extend $\overline{\varphi}_y$ to a Boolean homomorphism $\widehat{\varphi}_y$ onto $(\mathfrak{A} \widehat{\otimes}_R \mathfrak{B})/R$, the measure algebra of $\mathfrak{A} \widehat{\otimes}_R \mathfrak{B}$ modulo $R$ by [14], Theorem 33.1, page 141. Put $\varphi_y := \widehat{\varphi}_y \circ \widehat{\gamma}$ if $\widehat{\gamma} : \mathfrak{A} \widehat{\otimes}_R \mathfrak{B} \to (\mathfrak{A} \widehat{\otimes}_R \mathfrak{B})/R$ is the canonical surjection satisfying $\widehat{\gamma} | \mathfrak{A} \times \mathfrak{B} = \gamma$. Then $\varphi_y$ extends $\varphi_y^0$. □

THEOREM 2.2.  *Assume that a product r.c.p. $\{S_y : y \in Y\}$ on $\mathfrak{A}$ for $R$ with respect to $Q$ and a density $\rho \in \vartheta(Q)$ satisfy* (IT). *Then there exist $\psi \in \vartheta(\widehat{R})$ and $\psi_y \in \vartheta(\widehat{S}_y)$ such that the following conditions are satisfied for all $E \in \mathfrak{A} \widehat{\otimes}_R \mathfrak{B}$:*

 (i) $\widehat{S}_y([\psi(E)]^y \cup [\psi(E^c)]^y) = 1$ *for all $y \in Y$;*
 (ii) $\psi_y([\psi(E)]^y) = [\psi(E)]^y$ *for all $y \in Y$;*
 (iii) $[\psi(E)]_x \in \mathfrak{B}$ *for all $x \in X$;*
 (iv) $\psi(A \times B) \supseteq \bigcup_{y \in \rho(B)} (\psi_y(A) \times \{y\})$ *for all $A \in \mathfrak{A}$ and all $B \in \mathfrak{B}$.*

PROOF.  First recall that (IT) implies $S_y \ll P$ for all $y \in Y$. From $S_y \ll P$ for all $y \in Y$, we get, obviously, $R \ll P \otimes Q$ (see Section 1). If $f$ is a Radon–Nikodym derivative of $R$ with respect to $P \otimes Q$, put $E_R := \varphi(\{f > 0\})$. It follows that $E_R \in \mathfrak{A} \widehat{\otimes} \mathfrak{B}$ and $\varphi(E_R) = E_R$. Moreover, $(P \widehat{\otimes} Q)(E \cap E_R) = 0$ implies $\widehat{R}(E) = 0$ and $\widehat{R}(E_R) = 1$; that is, $E_R$ is the "measurable support" of the measure $R$. Moreover, $E_R^y \in \mathfrak{A}$ and $\sigma(E_R^y) = E_R^y$ for all $y \in Y$ by condition (b) above. For all $E \in \mathfrak{A} \widehat{\otimes} \mathfrak{B}$, define

$$\psi_1(E) := \begin{cases} X \times Y, & \text{if } E = X \times Y \text{ a.e. } (\widehat{R}), \\ \varphi(E \cap E_R), & \text{otherwise.} \end{cases}$$

By [8] we get $\psi_1 \in \vartheta(\widehat{R})$.

Next put $N := \{y \in Y : S_y(E_R^y) < 1\}$. Then $N \in \mathfrak{B}_0$. Notice that $E_R^y$ is a measurable support of $S_y$ for all $y \notin N$ since, for all such $y$, we have $S_y(E_R^y) = 1$ and $P(A \cap E_R^y) = 0$ implies $S_y(A) = 0$ because $S_y \ll P$. If, for all $A \in \mathfrak{A}$ and $y \notin N$, we define

$$\psi_y(A) := \begin{cases} X, & \text{if } A = X \text{ a.e. } (\widehat{S}_y), \\ \sigma(A \cap E_R^y), & \text{otherwise,} \end{cases}$$

then again $\psi_y \in \vartheta(S_y)$ for all $y \notin N$ by [8]. We denote the unique extension of $\psi_y$ to $\widehat{\mathfrak{A}}_y$ also by $\psi_y$.

For $y \in N$, take some $\tau_y \in \Lambda(\widehat{S}_y)$ and let $\varphi_y$ be defined according to Lemma 2.1. Then define $\psi_y := \tau_y$ for $y \in N$ and $\psi(E) := [\psi_1(E) \cap (X \times N^c)] \cup \bigcup_{y \in N} (\varphi_y(E) \times \{y\})$ for all $E \in \mathfrak{A} \widehat{\otimes}_R \mathfrak{B}$.

First $\psi(E) = \psi_1(E) = E$ a.e. $(\widehat{R})$ implies $\psi(E) = E$ a.e. $(\widehat{R})$ and $\psi(E) \in \mathfrak{A} \widehat{\otimes}_R \mathfrak{B}$. For $E, F \in \mathfrak{A} \widehat{\otimes}_R \mathfrak{B}$ with $E = F$ a.e. $(\widehat{R})$, we get $\psi_1(E \cap E_R) = \psi_1(F \cap E_R)$ and $\varphi_y(E) = \varphi_y(F)$, the latter by Lemma 2.1. This implies $\psi(E) = \psi(F)$. $\psi$ is stable with respect to intersections since $\psi_1$ and all the $\varphi_y$ are for all $y \in Y$. Since clearly $\psi(\varnothing) = \varnothing$, we get $\psi \in \vartheta(\widehat{R})$. So we are left with the verification of (i)–(iv).

(i) For $y \in N$, we have $\varphi_y(E) \cup \varphi_y(E^c) = X$ since $\varphi_y$ is a Boolean homomorphism: hence (i) is fullfilled in that case.

For $y \notin N$, we get $\widehat{S}_y([\psi(E)]^y \cup [\psi(E^c)]^y) = S_y(([\varphi(E)]^y \cup [\varphi(E^c)]^y) \cap E_R^y) = 1$, the latter since $P([\varphi(E)]^y \cup [\varphi(E^c)]^y) = 1$ for all $y \in Y$.

(ii) If $y \in N$, then $[\psi(E)]^y = \psi_y(A)$ for some $A \in \mathfrak{A}$, so $\psi_y([\psi(E)]^y) = \psi_y(\psi_y(A)) = \psi_y(A) = [\psi(E)]^y$.

If $y \in N^c$, then $[\psi(E)]^y = [\psi_1(E)]^y$. Hence if $\widehat{R}(X \times Y \triangle E) > 0$, then

$$\psi_y([\psi(E)]^y) = \psi_y([\psi_1(E)]^y) = \sigma([\psi_1(E)]^y \cap E_R^y)$$
$$= \sigma([\psi_1(E)]^y) \cap E_R^y = \sigma([\varphi(E)]^y \cap E_R^y) \cap E_R^y$$
$$= \sigma([\varphi(E)]^y) \cap E_R^y = [\varphi(E)]^y \cap E_R^y = [\varphi(E)]^y \cap [\varphi(E_R)]^y$$
$$= [\varphi(E \cap E_R)]^y = [\psi_1(E)]^y = [\psi(E)]^y.$$

Assertion (iii) immediately follows from the corresponding property of $\varphi$.

(iv) For $A \in \mathfrak{A}$ and all $B \in \mathfrak{B}$, we get $\psi(A \times B) = [((\sigma(A) \times \rho(B)) \cap E_R) \cap (X \times N^c)] \cup \bigcup_{y \in N}(\varphi_y(A \times B) \times \{y\})$.

For $y \in \rho(B) \cap N^c$, we get $[((\sigma(A) \times \rho(B)) \cap E_R) \cap (X \times N^c)]^y = \sigma(A) \cap E_R^y = \psi_y(A)$, that is, $[((\sigma(A) \times \rho(B)) \cap E_R) \cap (X \times N^c)] = \bigcup_{y \in \rho(B) \cap N^c}(\psi_y(A) \times \{y\})$.

For $y \in N$, we get $\varphi_y(A \times B) = \psi_y(A)$ if $y \in \rho(B)$ and $\varphi_y(A \times B) = \varnothing$ if $y \notin \rho(B)$. Both cases taken together give the assertion. □

THEOREM 2.3. *Assume that* (IT) *is satisfied by a product r.c.p.* $\{S_y : y \in Y\}$ *on* $\mathfrak{A}$ *for* $R$ *with respect to* $Q$ *and by* $\rho \in \Lambda(Q)$. *Then there exist* $\sigma_y \in \Lambda(\widehat{S}_y)$ *for all* $y \in Y$ *as well as* $\pi \in \Lambda(\widehat{R})$ *such that the following conditions hold true:*

(i) $[\pi(E)]^y = \sigma_y([\pi(E)]^y)$ *for all* $y \in Y$ *and* $E \in \mathfrak{A} \widehat{\otimes}_R \mathfrak{B}$;
(ii) $\pi(A \times B) = \bigcup_{y \in \rho(B)}(\sigma_y(A) \times \{y\})$ *for all* $A \in \mathfrak{A}$ *and all* $B \in \mathfrak{B}$.

PROOF. Let $\psi_y$ and $\psi$ be given as in Theorem 2.2. Next choose $\sigma_y \in \Lambda(\widehat{S}_y)$ satisfying $\psi_y(A) \subseteq \sigma_y(A)$ for all $A \in \mathfrak{A}$ and all $y \in Y$, and define $\pi \in \vartheta(\widehat{R})$ by setting, for each $E \in \mathfrak{A} \widehat{\otimes}_R \mathfrak{B}$ and each $y \in Y$,

$$[\pi(E)]^y = \sigma_y([\psi(E)]^y). \tag{3}$$

Since $\psi(E) \subseteq \pi(E)$ for all $E \in \mathfrak{A} \widehat{\otimes}_R \mathfrak{B}$, we get $\widehat{R}$-measurability of $\pi(E)$, and it follows easily that $\pi \in \vartheta(\widehat{R})$. In order to prove that $\pi$ is a lifting, it suffices to show that we have always $\pi(E^c) = [\pi(E)]^c$. But this is a consequence of Theorem 2.2 and (3), as we get for each $y$ the equality

$$[\pi(E^c)]^y = \sigma_y([\psi(E^c)]^y) = \sigma_y(([\psi(E)]^y)^c) = (\sigma_y([\psi(E)]^y))^c = ([\pi(E)]^y)^c.$$

This proves that $\pi \in \Lambda(\widehat{R})$. □

LEMMA 2.4. *If $R \ll P \otimes Q$, then for any $\rho \in \Lambda(Q)$ and for any product r.c.p. $\{S_y : y \in Y\}$ on $\mathfrak{A}$ for $R$ with respect to $Q$, there exist a product r.c.p. $\{T_y : y \in Y\}$ on $\mathfrak{A}$ for $R$ with respect to $Q$, $\varphi' \in \vartheta(\widehat{R})$ and $\rho' \in \Lambda(Q)$ such that $\{T_y \neq S_y\} \in \mathfrak{B}_0$ and (IT) is satisfied by the product r.c.p. $\{T_y : y \in Y\}$ and the lifting $\rho'$. Moreover, $\rho'$ and $\varphi'$ satisfy conditions (a)–(d) from the beginning of this section, too.*

PROOF. Put $N := \{y \in Y : S_y(E_R^y) < 1\}$ and note $N \in \mathfrak{B}_0$. Choose $y_0 \in N^c$ and define $\rho'(B) := [\rho(B) \cap N^c] \cup N$ if $y_0 \in \rho(B)$ and $\rho'(B) := \rho(B) \cap N^c$ otherwise. Put $\varphi'(E) := [\varphi(E) \cap (X \times N^c)] \cup \{(x,y) \in X \times N : (x, y_0) \in \varphi(E)\}$ for all $E \in \mathfrak{A} \widehat{\otimes}_R \mathfrak{B}$. It is now straightforward to verify that $\rho' \in \Lambda(Q)$.

To verify that $\varphi' \in \vartheta(\widehat{R})$, let $E, F \in \mathfrak{A} \widehat{\otimes}_R \mathfrak{B}$. The equality $\varphi'(E \cap F) = \varphi'(E) \cap \varphi'(F)$ is perhaps most easily verified by checking all equalities of sections $[\varphi'(E \cap F)]^y = [\varphi'(E)]^y \cap [\varphi'(F)]^y$ for all $y \in Y$. The other density properties of $\varphi'$ are straightforward to verify.

It is also easy to see that $\rho'$ and $\varphi'$ have all the properties (a)–(d) (with the same $\sigma$ as before). Next put $T_y := S_y$ for $y \in N^c$ and $T_y := S_{y_0}$ if $y \in N$. Clearly, $\{T_y : y \in Y\}$ is a product r.c.p. on $\mathfrak{A}$ for $R$ with respect to $Q$ such that $\{T_y \neq S_y\} \in \mathfrak{B}_0$. For $A \in \mathfrak{A}$ and $B \in \mathfrak{B}$, then $R(A \times B) = 0$ implies $\varphi'(A \times B) = \varnothing$, saying $\sigma(A) = \varnothing$ or $\rho'(B) = \varnothing$. But $\rho'(B) = \varnothing$ implies $Q(B) = 0$. If $Q(B) > 0$, then $\sigma(A) = \varnothing$, implying $T_y(A) = 0$ for all $y \in N^c$. But $T_y(A) = S_{y_0}(A) = 0$ for all $y \in N$ since $y_0 \in N^c$. So (IT) is satisfied for $\{T_y : y \in Y\}$ and $\rho'$. □

THEOREM 2.5. *If $R \ll P \otimes Q$, then there exist a product r.c.p. $\{T_y : y \in Y\}$ on $\mathfrak{A}$ for $R$ with respect to $Q$ and a lifting $\rho' \in \Lambda(Q)$ such that $\{T_y : y \in Y\}$ and $\rho'$ satisfy (IT), and there exist $\psi_y \in \vartheta(\widehat{T}_y)$ for all $y \in Y$ as well as $\psi \in \vartheta(\widehat{R})$ such that the following conditions hold true for all $E \in \mathfrak{A} \widehat{\otimes}_R \mathfrak{B}$:*

(i) $\widehat{T}_y([\psi(E)]^y \cup [\psi(E^c)]^y) = 1$ for all $y \in Y$;
(ii) $[\psi(E)]^y = \psi_y([\psi(E)]^y)$ for all $y \in Y$;
(iii) $[\psi(E)]_x \in \mathfrak{B}$ for all $x \in X$;
(iv) $\psi(A \times B) \supseteq \bigcup_{y \in \rho'(B)} (\psi_y(A) \times \{y\})$ for all $A \in \mathfrak{A}$ and all $B \in \mathfrak{B}$.

PROOF. We choose a product r.c.p. $\{S_y : y \in Y\}$ satisfying $S_y \ll P$ for all $y \in Y$. Now applying Lemma 2.4, we modify the product r.c.p. $\{S_y : y \in Y\}$, obtaining a product r.c.p. $\{T_y : y \in Y\}$, and find $\rho'$, $\varphi'$ satisfying all the conditions of Lemma 2.4. The product r.c.p. $\{T_y : y \in Y\}$ and the lifting $\rho'$ satisfy (IT). Now apply Theorem 2.2 for the product r.c.p. $\{T_y : y \in Y\}$ and for $\rho'$, instead of the product r.c.p. $\{S_y : y \in Y\}$ and $\rho$, and with $\varphi'$ instead of $\varphi$ in the proof. Theorem 2.2 produces now the required $\psi \in \vartheta(\widehat{R})$ and $\psi_y \in \vartheta(\widehat{T}_y), y \in Y$. □

THEOREM 2.6. *If $R \ll P \otimes Q$, then there exist a product r.c.p. $\{T_y : y \in Y\}$ on $\mathfrak{A}$ for $R$ with respect to $Q$ and a lifting $\rho' \in \Lambda(Q)$ such that $\{T_y : y \in Y\}$ and $\rho'$ satisfy* (IT), *and there exist $\sigma_y \in \Lambda(\widehat{T}_y)$ for all $y \in Y$ as well as $\varpi \in \Lambda(\widehat{R})$ such that the following conditions hold true:*

  (i) $[\varpi(E)]^y = \sigma_y([\varpi(E)]^y)$ *for all $y \in Y$ and all $E \in \mathfrak{A} \widehat{\otimes}_R \mathfrak{B}$;*
  (ii) $\varpi(A \times B) = \bigcup_{y \in \rho'(B)} (\sigma_y(A) \times \{y\})$ *for all $A \in \mathfrak{A}$ and all $B \in \mathfrak{B}$.*

PROOF. Let $T_y$, $\psi_y$, $\rho'$ and $\psi$ be given as in Theorem 2.5. Next choose $\sigma_y \in \Lambda(\widehat{T}_y)$ satisfying $\psi_y(A) \subseteq \sigma_y(A)$ for all $A \in \mathfrak{A}$ and all $y \in Y$ and define $\varpi \in \vartheta(\widehat{R})$ exactly as we have defined $\pi$ in the proof of Theorem 2.3. The required result follows in the same way as in the proof of Theorem 2.3. □

For given product r.c.p. $\{S_y : y \in Y\}$ on $\mathfrak{A}$ for $R$ with respect to $Q$, put $B_A := \{y \in Y : S_y(A) > 0\}$ for all $A \in \mathfrak{A}$. Define $\mathcal{B}_S := \{B_A : A \in \mathfrak{A}\}$. It follows that $Y \in \mathcal{B}_S$ and $\mathcal{B}_S$ is directed upwards since $B_A \cup B_C \subseteq B_{A \cup C}$ for $A, C \in \mathfrak{A}$. Denote by $\mathcal{T}_S$ the topology generated by $\mathcal{B}_S$ on $Y$ and let $\tau_\rho := \{B \in \mathfrak{B} : B \subseteq \rho(B)\}$ be one of the lifting topologies on $Y$ considered in [7], Chapter V. We do not assume in the next proposition that $R$ is absolutely continuous with respect to $P \otimes Q$.

PROPOSITION 2.7. *For a product r.c.p. $\{S_y : y \in Y\}$ on $\mathfrak{A}$ for $R$ with respect to $Q$, $\{\sigma_y \in \Lambda(\widehat{S}_y) : y \in Y\}$, $\rho \in \Lambda(Q)$ and $\pi \in \Lambda(\widehat{R})$, the following conditions are all equivalent:*

  (i) $\{S_y : y \in Y\}$ *and $\rho \in \Lambda(Q)$ satisfy the condition* (IT);
  (ii) $\mathcal{B}_S \subseteq \tau_\rho$;
  (iii) $\rho$ *is $\mathcal{T}_S$-strong;*
  (iv) $\rho, \{\sigma_y : y \in Y\}$ *and $\pi$ satisfy the rectangle formula* (RF).

PROOF. (i) ⇔ (ii) Assume that $B_A \subseteq \rho(B_A)$ for all $A \in \mathfrak{A}$ and take $A \in \mathfrak{A}$ and $B \in \mathfrak{B}$ such that $R(A \times B) = 0$. Then $\int_B S_y(A) \, dQ(y) = 0$. This yields $S_y(A) = 0$ for $Q$-almost all $y \in B$. Hence $\rho(B) \subseteq \rho(B_A^c)$. But by the assumption, we have $\rho(B_A^c) \subseteq B_A^c$ and so if $Q(B) > 0$, then $S_y(A) = 0$ provided $y \in \rho(B)$.

Now the reverse implication. Assume that (IT) is satisfied and take an $A \in \mathfrak{A}$. Then $R(A \times B_A^c) = 0$. If $Q(B_A^c) = 0$, then $Q(B_A) = 1$ and so $Y = \rho(B_A) \supseteq B_A$. If $Q(B_A^c) > 0$, then (IT) yields $S_y(A) = 0$ for all $y \in \rho(B_A^c)$. That means that $\rho(B_A^c) \subseteq B_A^c$, or equivalently, $B_A \subseteq \rho(B_A)$.

(ii) ⇔ (iii) is obvious since if $\mathcal{T}$ is an arbitrary topology on $Y$, then $\rho \in \Lambda(Q)$ is $\mathcal{T}$-strong for the topology $\mathcal{T}$ if and only if $\mathcal{T} \subseteq \tau_\rho$.

The implication (iv) ⇒ (i) is clear, while the implication (i) ⇒ (iv) follows from Theorem 2.3. □

**3. Arbitrary probability measures.** In this section we assume that we are given a lifting $\rho$ on $(Y, \mathfrak{B}, Q)$. We start with the following well-known result.

LEMMA 3.1. *Let $f$ be a bounded real-valued $(\mathfrak{A} \otimes \mathfrak{B})$-measurable function on $X \times Y$. Then the function $y \in Y \mapsto \int_X f^y(x)\, dS_y(x)$ is $\mathfrak{B}$-measurable and the equality*

$$\int_{X \times Y} f(x,y)\, dR(x,y) = \int_Y \int_X f^y(x)\, dS_y(x)\, dQ(y)$$

*holds true.*

*If $f$ is a bounded real-valued $(\mathfrak{A} \widehat{\otimes}_R \mathfrak{B})$-measurable function on $X \times Y$, then:*

(i) *the function $f^y$ is $\widehat{\mathfrak{A}}_y$-measurable for $Q$-a.a. $y \in Y$;*
(ii) *the function $y \mapsto \int_X f^y(x)\, d\widehat{S}_y(x)$ is $\widehat{\mathfrak{B}}$-measurable;*
(iii) *the equality*

$$\int_{X \times Y} f(x,y)\, d\widehat{R}(x,y) = \int_Y \int_X f^y(x)\, d\widehat{S}_y(x)\, d\widehat{Q}(y)$$

*holds true.*

*In particular, if $E \in \mathfrak{A} \widehat{\otimes}_R \mathfrak{B}$, then*

$$\widehat{R}(E) = \int_Y \widehat{S}_y(E^y)\, d\widehat{Q}(y).$$

The next result corresponds to Lemma 2.1 of [10] and plays an essential role in the proof of Proposition 3.4.

PROPOSITION 3.2. *If $\mathfrak{C}$ is a finite subalgebra of $\mathfrak{A}$, then for every $f \in \mathcal{L}^\infty(R) := \mathcal{L}^\infty(X \times Y, \mathfrak{A} \otimes \mathfrak{B}, R)$, we have the following section property of the conditional expectation:*

$$Y \setminus \{y \in Y : [E_{\mathfrak{C} \otimes \mathfrak{B}}(f)]^y = E_{\mathfrak{C}}^y(f^y) \quad a.e.\ (S_y|\mathfrak{C})\} \in \mathfrak{B}_0.$$

PROOF. By Lemma 3.1, if $f \in \mathcal{L}^\infty(R)$, $C \in \mathfrak{C}$, $B \in \mathfrak{B}$, then

$$\int_B \int_C [E_{\mathfrak{C} \otimes \mathfrak{B}}(f)]^y(x)\, dS_y(x)\, dQ(y)$$

$$= \int_{C \times B} E_{\mathfrak{C} \otimes \mathfrak{B}}(f)(x,y)\, dR(x,y)$$

$$= \int_{C \times B} f(x,y)\, dR(x,y)$$

$$= \int_B \int_C f^y(x)\, dS_y(x)\, dQ(y)$$

$$= \int_B \int_C E_{\mathfrak{C}}^y(f^y)(x)\,dS_y(x)\,dQ(y).$$

This implies that, for each $C \in \mathfrak{C}$, there exists $N_C \in \mathfrak{B}_0$ such that, for any $y \in Y \setminus N_C$ we have

$$\int_C [E_{\mathfrak{C} \otimes \mathfrak{B}}(f)]^y(x)\,dS_y(x) = \int_C E_{\mathfrak{C}}^y(f^y)(x)\,dS_y(x).$$

If $N_f := \bigcup_{C \in \mathfrak{C}} N_C$, then $N_f \in \mathfrak{B}_0$, and for all $y \in Y \setminus N_f$ and all $C \in \mathfrak{C}$, we have

$$\int_C [E_{\mathfrak{C} \otimes \mathfrak{B}}(f)]^y(x)\,dS_y(x) = \int_C E_{\mathfrak{C}}^y(f^y)(x)\,dS_y(x).$$

It follows that we have, for all $y \in Y \setminus N_f$,

$$[E_{\mathfrak{C} \otimes \mathfrak{B}}(f)]^y = E_{\mathfrak{C}}^y(f^y) \qquad \text{a.e. } (S_y|\mathfrak{C}). \qquad \square$$

LEMMA 3.3. *Let $\mathfrak{C}$ be a finite subalgebra of $\mathfrak{A}$ and let $H \in \mathfrak{A} \setminus \mathfrak{C}$. Let $\mathfrak{D} = \sigma(\mathfrak{C} \cup \{H\})$. Assume that for each $y \in Y$, there exists $\tau_y \in \Lambda(S_y|\mathfrak{C})$ and $\varphi \in \Lambda(R|\mathfrak{C} \otimes \mathfrak{B})$ such that*

(4) $[\varphi(A \times B)]^y = \tau_y(A) \qquad$ *for all $A \in \mathfrak{C}, B \in \mathfrak{B}$ and $Q$-almost all $y \in B$,*

(5) $[\varphi(A \times B)]^y = \varnothing \qquad$ *for all $A \in \mathfrak{C}, B \in \mathfrak{B}$ and $Q$-almost all $y \in B^c$.*

*Then there exist $\xi \in \Lambda(R|\mathfrak{D} \otimes \mathfrak{B})$ and $\xi_y \in \Lambda(S_y|\mathfrak{D})$ for $y \in Y$ such that*

(6) $[\xi(A \times B)]^y = \xi_y(A) \qquad$ *for all $A \in \mathfrak{D}, B \in \mathfrak{B}$ and $Q$-almost all $y \in B$,*

(7) $[\xi(A \times B)]^y = \varnothing \qquad$ *for all $A \in \mathfrak{D}, B \in \mathfrak{B}$ and $Q$-almost all $y \in B^c$.*

*Moreover, $\xi|\mathfrak{C} \otimes \mathfrak{B} = \varphi$ and $\xi_y|\mathfrak{C} = \tau_y$ for all $y \in Y$.*

PROOF. By induction, one can always assume that $H$ is contained in some atom (call it $A_H$) of $\mathfrak{C}$. We set then

$$\xi_y(H) = \begin{cases} \varnothing, & \text{if } S_y(H) = 0, \\ \tau_y(A_H), & \text{if } S_y(A_H \setminus H) = 0, \\ H, & \text{otherwise,} \end{cases}$$

and

$$\xi_y(A_H \setminus H) = \begin{cases} \tau_y(A_H), & \text{if } S_y(H) = 0, \\ \varnothing, & \text{if } S_y(A_H \setminus H) = 0, \\ \tau_y(A_H) \setminus H, & \text{otherwise.} \end{cases}$$

If $D = (A \cap H) \cup (B \cap H^c)$ with $A, B \in \mathfrak{C}$ given, then we set

$$\xi_y(D) = (\xi_y(H) \cap \tau_y(A)) \cup (\xi_y(A_H \setminus H) \cap \tau_y(B)) \cup \tau_y(B \cap A_H^c).$$

We may write it also in a more symmetric way, if we notice that $\xi_y(H^c) = \xi_y(A_H \setminus H) \cup \xi_y(A_H^c)$. Then we have

$$\xi_y(D) = (\xi_y(H) \cap \tau_y(A)) \cup (\xi_y(H^c) \cap \tau_y(B)).$$

Similarly, if $B_A := \{y \in Y : S_y(A) > 0\}$, then we set

$$\xi(H \times Y) := [(H \times Y) \cap \varphi(A_H \times B_H)]$$
$$\cup [(H^c \times Y) \cap \varphi(A_H \times (Y \setminus B_{A_H \setminus H}))],$$
$$\xi((A_H \setminus H) \times Y) := [(H \times Y) \cap \varphi(A_H \times B_H^c)]$$
$$\cup [(H^c \times Y) \cap \varphi(A_H \times B_{A_H \setminus H})],$$

and if $E = [F \cap (H \times Y)] \cup G \cap [(H^c \times Y)]$ with $F, G \in \mathfrak{C} \otimes \mathfrak{B}$, then

$$\xi(E) = (\varphi(F) \cap \xi(H \times Y)) \cup (\varphi(G) \cap \xi((A_H \setminus H) \times Y)) \cup \varphi((G \cap (A_H^c \times Y))$$

or

$$\xi(E) = (\varphi(F) \cap \xi(H \times Y)) \cup (\varphi(G) \cap \xi(H^c \times Y)).$$

We are now going to check the section properties of $\xi$. If $y \in Y$, then

$$[\xi(H \times Y)]^y = (H \cap [\varphi(A_H \times B_H)]^y) \cup (H^c \cap [\varphi(A_H \times (Y \setminus B_{A_H \setminus H}))]^y).$$

Now

$$[\varphi(A_H \times B_H)]^y = \begin{cases} \tau_y(A_H), & \text{for } Q\text{-a.e. } y \in B_H, \\ \varnothing, & \text{for } Q\text{-a.e.} y \notin B_H, \end{cases}$$

and

$$[\varphi(A_H \times (Y \setminus B_{A_H \setminus H}))]^y = \begin{cases} \tau_y(A_H), & \text{for } Q\text{-a.e. } y \notin B_{A_H \setminus H}, \\ \varnothing, & \text{for } Q\text{-a.e. } y \in B_{A_H \setminus H}. \end{cases}$$

The above relations give $[\xi(H \times Y)]^y = \xi_y(H)$ for $Q$-almost all $y \in Y$.

In a similar way, we get

$$[\xi(H^c \times Y)]^y = \xi_y(H^c)$$

for $Q$-almost all $y \in Y$.

Consider now arbitrary $A \in \mathfrak{D}$, $B \in \mathfrak{B}$. Then there exist $\widetilde{A}_1, \widetilde{A}_2 \in \mathfrak{C}$ such that

$$A = (\widetilde{A}_1 \cap H) \cup (\widetilde{A}_2 \cap H^c)$$

and so

$$A \times B = [(\widetilde{A}_1 \times B) \cap (H \times Y)] \cup [(\widetilde{A}_2 \times B) \cap (H^c \times Y)].$$

Then, applying (4) and (5), we have, for $Q$-almost all $y \in B$,

$$\begin{aligned}[\xi(A \times B)]^y &= ([\varphi(\widetilde{A}_1 \times B)]^y \cap [\xi(H \times Y)]^y) \\ &\quad \cup ([\varphi(\widetilde{A}_2 \times B)]^y \cap [\xi(H^c \times Y)]^y) \\ &= ([\varphi(\widetilde{A}_1 \times B)]^y \cap \xi_y(H)) \cup ([\varphi(\widetilde{A}_2 \times B)]^y \cap \xi_y(H^c)) \\ &= ([\tau_y(\widetilde{A}_1)]^y \cap \xi_y(H)) \cup ([\tau_y(\widetilde{A}_2)]^y \cap \xi_y(H^c)) \\ &= \xi_y(A)\end{aligned}$$

and $[\xi(A \times B)]^y = \varnothing$ for $Q$-almost all $y \in B^c$. That proves (6) and (7). $\square$

PROPOSITION 3.4. *Assume that $\mathfrak{A}$ contains a countably generated $\sigma$-algebra which is dense in $\mathfrak{A}$ (in the Fréchet–Nikodym pseudometric) with respect to $P$. Then there exist $\widetilde{\varphi} \in \vartheta(R)$ and $\{\tau_y \in \vartheta(S_y) : y \in Y\}$ such that, for each $F \in \mathfrak{A} \otimes \mathfrak{B}$,*

$$[\widetilde{\varphi}(F)]^y = \tau_y([\widetilde{\varphi}(F)]^y) \qquad \textit{for almost all } y \in Y.$$

*There exists also $\varphi \in \vartheta(\widehat{R})$ such that, for each $F \in \mathfrak{A} \widehat{\otimes}_R \mathfrak{B}$,*

$$[\varphi(F)]^y = \tau_y([\varphi(F)]^y) \qquad \textit{for all } y \in Y$$

*and*

$$[\varphi(F)]_x \in \mathfrak{B} \qquad \textit{for all } x \in X.$$

PROOF. Let $(\mathfrak{C}_n)_{n \in \mathbf{N}}$ be a sequence of finite algebras $\mathfrak{C}_n$ such that $\mathfrak{C}_1 = \{\varnothing, X\}$, $\mathfrak{C}_n \subseteq \mathfrak{C}_{n+1}$ for $n \in \mathbf{N}$ and $\mathfrak{D} := \sigma(\bigcup_{n \in \mathbf{N}} \mathfrak{C}_n)$ is $P$-dense in $\mathfrak{A}$. It is easily seen that $\mathfrak{D} \otimes \mathfrak{B} = \sigma(\bigcup_{n \in \mathbf{N}} (\mathfrak{C}_n \otimes \mathfrak{B}))$.

Then, we construct $\varphi_n \in \Lambda(R | \mathfrak{C}_n \otimes \mathfrak{B})$ and $\tau_{ny} \in \Lambda(S_y | \mathfrak{C}_n)$ according to Lemma 3.3, taking as $\varphi_1 \in \Lambda(R | \mathfrak{C}_1 \otimes \mathfrak{B})$ the lifting determined by $\rho \in \Lambda(Q)$ and by $\tau_{1y} \in \Lambda(S_y | \mathfrak{C}_1)$ the only existing lifting on $(X, \mathfrak{C}_1, S_y | \mathfrak{C}_1)$. Following [6], we define $\widetilde{\varphi} \in \vartheta(R | \mathfrak{D} \otimes \mathfrak{B})$ and $\widetilde{\tau}_y \in \vartheta(S_y | \mathfrak{D})$ with $\widetilde{\varphi} | (\mathfrak{C}_n \otimes \mathfrak{B}) = \varphi_n$ and $\widetilde{\tau}_y | \mathfrak{C}_n = \tau_{ny}$ for all $y \in Y$ and all $n \in \mathbf{N}$ by means of

$$\widetilde{\varphi}(F) := \bigcap_{k \in \mathbf{N}} \bigcup_{n \in \mathbf{N}} \bigcap_{m \geq n} \varphi_m(\{E_{\mathfrak{C}_m \otimes \mathfrak{B}}(\chi_F) > 1 - 1/k\})$$

if $F \in \mathfrak{D} \otimes \mathfrak{B}$ and

$$\widetilde{\tau}_y(A) := \bigcap_{k \in \mathbf{N}} \bigcup_{n \in \mathbf{N}} \bigcap_{m \geq n} \tau_{my}(\{E^y_{\mathfrak{C}_m}(\chi_A) > 1 - 1/k\})$$

if $A \in \mathfrak{D}$. Then, for $F \in \mathfrak{D} \otimes \mathfrak{B}$, let $N_{m,F} \in \mathfrak{B}_0$ be such that (see Proposition 3.2)

$$\{y \in Y : [E_{\mathfrak{C}_m \otimes \mathfrak{B}}(\chi_F)]^y = E^y_{\mathfrak{C}_m}(\chi_{F^y}) \text{ a.e.}(S_y | \mathfrak{C}_m)\}^c \subseteq N_{m,F}.$$

If $N_F := \bigcup_{m \in \mathbf{N}} N_{m,F}$, then $N_F \in \mathfrak{B}_0$ and it follows from Proposition 3.2 that, for all $y \in Y \setminus N_F$,

$$[\widetilde{\varphi}(F)]^y = \bigcap_{k \in \mathbf{N}} \bigcup_{n \in \mathbf{N}} \bigcap_{m \geq n} [\varphi_m(\{E_{\mathfrak{C}_m \otimes \mathfrak{B}}(\chi_F) > 1 - 1/k\})]^y$$

$$= \bigcap_{k \in \mathbf{N}} \bigcup_{n \in \mathbf{N}} \bigcap_{m \geq n} \tau_{my}([\{E_{\mathfrak{C}_m \otimes \mathfrak{B}}(\chi_F) > 1 - 1/k\}]^y)$$

$$= \bigcap_{k \in \mathbf{N}} \bigcup_{n \in \mathbf{N}} \bigcap_{m \geq n} \tau_{my}(\{[E_{\mathfrak{C}_m \otimes \mathfrak{B}}(\chi_F)]^y > 1 - 1/k\})$$

$$= \bigcap_{k \in \mathbf{N}} \bigcup_{n \in \mathbf{N}} \bigcap_{m \geq n} \tau_{my}(\{E^y_{\mathfrak{C}_m}(\chi_{F^y}) > 1 - 1/k\}) = \widetilde{\tau}_y(F^y).$$

It follows that $\widetilde{\tau}_y([\widetilde{\varphi}(F)]^y) = [\widetilde{\varphi}(F)]^y$ for all $y \notin N_F$.

Notice that according to Lemma 3.3, the sets $\varphi_m(\{E_{\mathfrak{C}_m \otimes \mathfrak{B}}(\chi_F) > 1 - 1/k\})$ are members of $\mathfrak{C}_m \otimes \mathfrak{B}$ and so the sets $\tau_{my}([\{E_{\mathfrak{C}_m \otimes \mathfrak{B}}(\chi_F) > 1 - 1/k\}]^y)$ are properly defined. Then, we set $[\varphi(F)]^y = \widetilde{\tau}_y([\widetilde{\varphi}(F)]^y)$ for every $y \in Y$ and $F \in \mathfrak{D} \otimes \mathfrak{B}$. Notice that since $\widetilde{\varphi}(F) \in \mathfrak{D} \otimes \mathfrak{B}$, we have $[\widetilde{\varphi}(F)]^y \in \mathfrak{D}$ and $[\widetilde{\varphi}(F)]_x \in \mathfrak{B}$, for every $y \in Y$ and $x \in X$. This proves the correctness of the definition of $\varphi$. It is clear that $\varphi(F) \in \mathfrak{A} \widehat{\otimes}_R \mathfrak{B}$. But as $\mathfrak{D} \otimes \mathfrak{B}$ is $R$-dense in $\mathfrak{A} \otimes \mathfrak{B}$, we can extend $\varphi$ in the unique way to the whole $\mathfrak{A} \widehat{\otimes}_R \mathfrak{B}$. We denote this extension also by $\varphi$.

At the moment, each density $\widetilde{\tau}_y$ is defined on the $\sigma$-algebra $\mathfrak{D}$. For every $y \in Y$, denote again by $\widetilde{\tau}_y$ the obvious extension of $\widetilde{\tau}_y$ to $\mathfrak{D}_y := \sigma(\mathfrak{D} \cup \mathfrak{A}_{y0})$. It is well known (cf. [6]) that each density on a $\sigma$-subalgebra containing all sets of measure zero can be extended to a density on the whole $\sigma$-algebra without changing its values on the smaller $\sigma$-algebra. Let $\tau_y \in \vartheta(S_y)$ be an arbitrary extension of $\widetilde{\tau}_y$ from $\mathfrak{D}_y$ to $\mathfrak{A}$. One can easily see that $\varphi$ and $\{\tau_y : y \in Y\}$ satisfy the required properties. $\square$

THEOREM 3.5. *Assume that $\mathfrak{A}$ contains a countably generated $\sigma$-algebra which is dense in $\mathfrak{A}$ (in the Fréchet–Nikodym pseudometric) with respect to $P$. Then there exist $\tau_y \in \vartheta(\widehat{S}_y)$ for all $y \in Y$ and $\psi \in \vartheta(\widehat{R})$ satisfying, for all $E \in \mathfrak{A} \widehat{\otimes}_R \mathfrak{B}$, the following conditions:*

(i) $\widehat{S}_y([\psi(E)]^y \cup [\psi(E^c)]^y) = 1$ *for all $y \in Y$;*
(ii) $[\psi(E)]^y = \tau_y([\psi(E)]^y)$ *for all $y \in Y$;*
(iii) $[\psi(E)]_x$ *is $\widehat{Q}$-measurable for all $x \in X$.*

PROOF. There exist $\varphi \in \vartheta(\widehat{R})$ and $\tau_y \in \vartheta(S_y)$ for all $y \in Y$ satisfying the thesis of Proposition 3.4. We denote again by $\tau_y$ the unique extension of $\tau_y$ to $\widehat{\mathfrak{A}}_y$. Let

$$\Phi := \{\overline{\varphi} \in \vartheta(\widehat{R}) : \forall y \in Y \ \forall E \in \mathfrak{A} \widehat{\otimes}_R \mathfrak{B} \ [\overline{\varphi}(E)]^y \subseteq \tau_y([\overline{\varphi}(E)]^y)$$
$$\text{and} \ \forall E \in \mathfrak{A} \widehat{\otimes}_R \mathfrak{B} \ \varphi(E) \subseteq \overline{\varphi}(E)\}.$$

Notice first that $\Phi \neq \varnothing$ since $\varphi \in \Phi$.

We consider $\Phi$ with inclusion as the partial order: $\overline{\varphi}_1 \leq \overline{\varphi}_2$ if $\overline{\varphi}_1(E) \subseteq \overline{\varphi}_2(E)$ for each $E \in \mathfrak{A}\widehat{\otimes}_R \mathfrak{B}$.

CLAIM 1. *There exists a maximal element in $\Phi$.*

PROOF. The only fact we have to prove is that each chain $\{\varphi_\alpha\}_{\alpha \in A} \subseteq \Phi$ has a dominating element in $\Phi$. The obvious candidate is $\overline{\varphi}$ given, for each $E \in \mathfrak{A}\widehat{\otimes}_R \mathfrak{B}$, by

$$\overline{\varphi}(E) = \bigcup_{\alpha \in A} \varphi_\alpha(E). \tag{8}$$

It can be easily seen that $\overline{\varphi}$ is a density dominating the chain and $\overline{\varphi} \in \Phi$.

Indeed, let us prove first the measurability of $\overline{\varphi}(E)$. To do it notice first that

$$\overline{\varphi}(E^c) = \bigcup_{\alpha \in A} \varphi_\alpha(E^c). \tag{9}$$

and that (8) and (9) together yield immediately

$$\overline{\varphi}(E) \cap \overline{\varphi}(E^c) = \varnothing.$$

Hence, $\overline{\varphi}(E) \subseteq [\overline{\varphi}(E^c)]^c$ and so if an $\alpha \in A$ is fixed, then

$$\varphi_\alpha(E) \subseteq \overline{\varphi}(E) \subseteq [\overline{\varphi}(E^c)]^c \subseteq [\varphi_\alpha(E^c)]^c$$

for each $E \in \mathfrak{A}\widehat{\otimes}_R \mathfrak{B}$. Since $\widehat{R}$ is complete and $\varphi_\alpha \in \vartheta(\widehat{R})$, this proves the measurability of $\overline{\varphi}(E)$ and the relation $\overline{\varphi}(E) \stackrel{\widehat{R}}{=} E$. Consider now the section properties of $\overline{\varphi}(E)$. For fixed $y \in Y$,

$$[\overline{\varphi}(E)]^y = \bigcup_{\alpha \in A} [\varphi_\alpha(E)]^y \subseteq \bigcup_{\alpha \in A} \tau_y([\varphi_\alpha(E)]^y),$$

and so, by virtue of [7], Chapter III, Section 3, the set $[\overline{\varphi}(E)]^y$ is $\widehat{S}_y$-measurable. Setting in the above inclusions $\overline{\varphi}(E)$ instead of $\varphi_\alpha(E)$, we see that the inclusion $[\overline{\varphi}(E)]^y \subseteq \tau_y([\overline{\varphi}(E)]^y)$ is satisfied also.

We have to prove yet that $\overline{\varphi}$ is a density. To do it, take $E, F \in \mathfrak{A}\widehat{\otimes}_R \mathfrak{B}$. We have

$$\overline{\varphi}(E) \cap \overline{\varphi}(F) = \bigcup_{\alpha \in A} \varphi_\alpha(E) \cap \bigcup_{\alpha \in A} \varphi_\alpha(F)$$

$$= \bigcup_{\alpha \in A} \left[ \varphi_\alpha(E) \cap \bigcup_{\beta \in A} \varphi_\beta(F) \right]$$

$$= \bigcup_{\alpha \in A} \left[ \bigcup_{\beta \in A} \varphi_\alpha(E) \cap \varphi_\beta(F) \right]$$

$$= \bigcup_{\alpha \in A} \left[ \bigcup_{\alpha \leq \beta} \varphi_\alpha(E) \cap \varphi_\beta(F) \right] \subseteq \bigcup_{\alpha \in A} \left[ \bigcup_{\alpha \leq \beta} \varphi_\beta(E) \cap \varphi_\beta(F) \right]$$

$$= \bigcup_{\alpha \in A} \left[ \bigcup_{\alpha \leq \beta} \varphi_\beta(E \cap F) \right] \subseteq \overline{\varphi}(E \cap F).$$

The reverse inclusion and other properties are clear and so $\overline{\varphi} \in \vartheta(\widehat{R})$. This proves that $\overline{\varphi}$ dominates the whole chain. According to the Zorn–Kuratowski lemma the set $\Phi$ possesses a maximal element $\psi$. $\square$

CLAIM 2. *For each $y \in Y$ and $E \in \mathfrak{A} \widehat{\otimes}_R \mathfrak{B}$,*

$$\widehat{S}_y([\psi(E)]^y \cup [\psi(E^c)]^y) = 1.$$

PROOF. Notice first that as $\psi \in \Phi$, all sections $[\psi(E)]^y$ are $\widehat{S}_y$-measurable. Suppose now that there exist $H \in \mathfrak{A} \widehat{\otimes}_R \mathfrak{B}$ and $y_0 \in Y$ such that $\widehat{S}_{y_0}([\psi(H)]^{y_0} \cup [\psi(H^c)]^{y_0}) < 1$. Let

$$W := \tau_{y_0}[([\psi(H)]^{y_0} \cup [\psi(H^c)]^{y_0})^c]$$

and let

$$[\widehat{\psi}(E)]^y = \begin{cases} [\psi(E)]^y, & \text{if } y \neq y_0, \\ [\psi(E)]^{y_0} \cup (W \cap [\psi(H \cup E)]^{y_0}), & \text{if } y = y_0. \end{cases}$$

It is clear that $\psi(E) \subseteq \widehat{\psi}(E)$ for each $E \in \mathfrak{A} \widehat{\otimes}_R \mathfrak{B}$. In particular, $\widehat{\psi}(X \times Y) = X \times Y$. Also the other density properties are fulfilled.

It follows directly from the definition that $\widehat{\psi}$ and $\psi$ are different densities. In order to get a contradiction with our hypothesis, it is enough to show that $[\widehat{\psi}(E)]^{y_0} \subseteq \tau_{y_0}([\widehat{\psi}(E)]^{y_0})$, but this is immediate. If $E \in \mathfrak{A} \widehat{\otimes}_R \mathfrak{B}$, then

$$\tau_{y_0}([\widehat{\psi}(E)]^{y_0}) \supseteq \tau_{y_0}([\psi(E)]^{y_0}) \cup \tau_{y_0}(W \cap [\psi(H \cup E)]^{y_0})$$
$$\supseteq [\psi(E)]^{y_0} \cup [\tau_{y_0}(W) \cap ([\psi(H \cup E)]^{y_0})]$$
$$\supseteq [\psi(E)]^{y_0} \cup (W \cap [\psi(H \cup E)]^{y_0})$$
$$= [\widehat{\psi}(E)]^{y_0}.$$

This completes the proof of the claim. $\square$

CLAIM 3. *For each $y \in Y$ and $E \in \mathfrak{A} \widehat{\otimes}_R \mathfrak{B}$,*

$$[\psi(E)]^y = \tau_y([\psi(E)]^y).$$

PROOF. Set, for each $y \in Y$ and $E \in \mathfrak{A} \widehat{\otimes}_R \mathfrak{B}$,
$$[\widetilde{\psi}(E)]^y = \tau_y([\psi(E)]^y).$$
Clearly, $\psi(F) \subseteq \widetilde{\psi}(F)$ for each $F$. Moreover, the equality $\psi(E) \cap \psi(E^c) = \varnothing$ yields for each $y$ the relation $\tau_y([\psi(E)]^y) \cap \tau_y([\psi(E^c)]^y) = \varnothing$. As a consequence, we get $\widetilde{\psi}(E) \cap \widetilde{\psi}(E^c) = \varnothing$, and then $\widetilde{\psi}(E^c) \subseteq (\widetilde{\psi}(E))^c$. Hence
$$\psi(E^c) \subseteq \widetilde{\psi}(E^c) \subseteq [\widetilde{\psi}(E)]^c \subseteq [\psi(E)]^c.$$
Since $\psi \in \vartheta(\widehat{R})$, we have $\widehat{R}([\psi(E)]^c) = \widehat{R}[\psi(E^c)]$ and so $\widetilde{\psi}(E)$ is $\widehat{R}$-measurable. It follows that $\widetilde{\psi} \in \Phi$ and so $\psi = \widetilde{\psi}$, due to the maximality of $\psi$. □

CLAIM 4. *$X$-sections of $\psi$ are $\widehat{Q}$-measurable.*

PROOF. In order to prove the measurability of the $X$-sections of $\psi$, notice that since $\varphi$ and $\psi$ are densities in the same measure space, the equalities $\widehat{R}(\varphi(E) \triangle \psi(E)) = 0$ and $\widehat{R}(\varphi(E) \cup \varphi(E^c)) = 1$ hold true for every $E \in \mathfrak{A} \widehat{\otimes}_R \mathfrak{B}$. It follows then from Lemma 3.1 that there is $M_E \in \mathfrak{B}_0$ such that, for all $y \notin M_E$,
$$\widehat{S}_y([\varphi(E)]^y \triangle [\psi(E)]^y) = 0 \quad \text{and} \quad \widehat{S}_y([\varphi(E)]^y \cup [\varphi(E^c)]^y) = 1.$$
It follows that if $y \notin M_E$, then
$$[\psi(E)]^y = \tau_y([\psi(E)]^y) = \tau_y([\varphi(E)]^y) = [\varphi(E)]^y.$$
Hence,
$$\varphi(E) \cap (X \times M_E^c) = \psi(E) \cap (X \times M_E^c),$$
and consequently, we have for all $x \in X$,
$$[\varphi(E)]_x \cap M_E^c = [\psi(E)]_x \cap M_E^c.$$
Since all sections $[\varphi(E)]_x$ are $\widehat{Q}$-measurable, it follows that the sections $[\psi(E)]_x$ are $\widehat{Q}$-measurable. □

This completes the proof of Theorem 3.5. □

THEOREM 3.6. *Assume that $\mathfrak{A}$ contains a countably generated $\sigma$-algebra which is dense in $\mathfrak{A}$ (in the Fréchet–Nikodym pseudometric) with respect to $P$. Then there exist $\sigma_y \in \Lambda(\widehat{S}_y)$ for all $y \in Y$ and $\pi \in \Lambda(\widehat{R})$ such that the following condition is satisfied:*

(10) $\quad [\pi(E)]^y = \sigma_y([\pi(E)]^y) \quad$ *for all $y \in Y$ and $E \in \mathfrak{A} \widehat{\otimes}_R \mathfrak{B}$.*

*Equivalently, for each $f \in \mathcal{L}^\infty(\widehat{R})$ and each $y \in Y$,*
$$[\pi(f)]^y = \sigma_y([\pi(f)]^y).$$

PROOF. According to Theorem 3.5, there exist $\tau_y \in \vartheta(\widehat{S}_y)$ for all $y \in Y$ and $\psi \in \vartheta(\widehat{R})$ such that, for all $E \in \mathfrak{A}\widehat{\otimes}_R \mathfrak{B}$,

$$[\psi(E)]^y = \tau_y([\psi(E)]^y) \qquad \text{for all } y \in Y \tag{11}$$

and

$$\widehat{S}_y([\psi(E)]^y \cup [\psi(E^c)]^y) = 1 \qquad \text{for all } y \in Y. \tag{12}$$

We take now, for each $y \in Y$, a lifting $\sigma_y \in \Lambda(\widehat{S}_y)$ such that $\tau_y \subseteq \sigma_y$ and define $\pi \in \vartheta(\widehat{R})$ by setting, for each $E \in \mathfrak{A}\widehat{\otimes}_R \mathfrak{B}$ and each $y \in Y$,

$$[\pi(E)]^y = \sigma_y([\psi(E)]^y). \tag{13}$$

Since $\psi(E) \subseteq \pi(E)$ for all $E \in \mathfrak{A}\widehat{\otimes}_R \mathfrak{B}$, we get $\widehat{R}$-measurability of $\pi(E)$ and $\pi(E) \stackrel{\widehat{R}}{=} E$. In order to prove that $\pi$ is a lifting, it suffices to show that we have always $\pi(E^c) = [\pi(E)]^c$. But this is a consequence of (12) and (13) as we get for each $y$ the equality $[\pi(E^c)]^y = ([\pi(E)]^y)^c$. This proves that $\pi \in \Lambda(\widehat{R})$. □

REMARK 3.7. There is an obvious question: Can lifting $\pi$ (or $\varpi$ from Theorem 2.6) be chosen in such a way that all the sections $[\pi(E)]_x$ would be $\widehat{Q}$-measurable? Such a property holds true for the density $\psi$ in Theorem 3.5. As it has been observed in [10], an improvement of this type is in general impossible in case of product measures $R$.

**4. Examples.** In case of the regular conditional probabilities defined on the same basic space where $\mathfrak{B} \subset \mathfrak{A}$, the following result holds true (see [2] or [9], page 358):

Let $(X, \mathfrak{A}, P)$ be a probability space and let $\mathfrak{B}$ be a sub-$\sigma$-algebra of $\mathfrak{A}$. Assume that $\{P_x : x \in X\}$ is a r.c.p. on $\mathfrak{A}$ with respect to $\mathfrak{B}$. Then, if there exists a probability measure $\widetilde{P}$ on $\mathfrak{A}$ such that every measure $P_x$ is absolutely continuous with respect to $\widetilde{P}$, then $P|\mathfrak{B}$ is atomic and, for each $A \in \mathfrak{A}$, there exists $N_A \in \mathfrak{B}_0$ such that

$$P_x(A) = \sum_n \frac{P(A \cap B_n)}{P(B_n)} \chi_{B_n}(x) \qquad \text{for every } x \in X \setminus N_A,$$

where $B_1, B_2, \ldots$ are all the atoms of $P|\mathfrak{B}$.

If $\mathfrak{A}$ is countably generated, then one can replace the sets $N_A$ by one set $N \in \mathfrak{B}_0$ satisfying the above equality for all $A \in \mathfrak{A}$.

One may ask if a similar simplification takes place also in case of our investigations. The following example shows that this is not the case. When the algebras $\mathfrak{A}$ and $\mathfrak{B}$ are "independent," then there are nice examples with absolutely continuous product regular conditional probabilities. The product r.c.p. in the first example is even uniformly absolutely continuous, that is,

$$\forall \varepsilon > 0 \ \exists \delta > 0 \ \forall A \in \mathfrak{A}, \qquad \left[ P(A) < \delta \ \implies \ \sup_{y \in Y} S_y(A) < \varepsilon \right].$$

EXAMPLE 4.1. Let $(X, \mathfrak{A}) = (Y, \mathfrak{B}) = ((-\infty, +\infty), \mathfrak{L})$, where $\mathfrak{L}$ denotes the Lebesgue measurable sets, and let $\lambda$ be the Lebesgue measure. Moreover, let $C \subset [0, 1]$ be a set of positive measure: $0 < \lambda(C) \leq 1$. For each $y \in Y$, let $S_y$ be the measure defined on $\mathfrak{A}$ by

$$S_y(A) := \frac{1/\sqrt{2\pi} \int_{A \setminus (C+y)} \exp\{-(t-y)^2/2\} \, dt}{1/\sqrt{2\pi} \int_{R \setminus (C+y)} \exp\{-(t-y)^2/2\} \, dt}$$

$$= \frac{\int_{A \setminus (C+y)} \exp\{-(t-y)^2/2\} \, dt}{\int_{R \setminus C} \exp\{-t^2/2\} \, dt}$$

if $A \in \mathfrak{A}$. Then set

$$Q(B) = \frac{1}{\sqrt{2\pi}} \int_B \exp(-t^2/2) \, dt,$$

$$R(A \times B) = \int_B S_y(A) \, dQ(y)$$

and

$$P(A) = R(A \times Y).$$

$R$ is obviously additive on $\mathfrak{A} \times \mathfrak{B}$, and since $P$ and $Q$ are perfect, $R$ is countably additive on $\mathfrak{A} \times \mathfrak{B}$ and so it can be uniquely extended to a measure on $\mathfrak{A} \otimes \mathfrak{B}$ (see [13]). We denote the extension also by $R$.

CLAIM 1. *For every bounded $A \in \mathfrak{A}$, the function $S.(A)$ is continuous and $Q \ll P$.*

CLAIM 2. *If $\rho$ is a strong lifting for $Q$, then $\{S_y : y \in Y\}$ satisfies* (IT).

PROOF. If $R(A \times B) = 0$ and $Q(B) > 0$, then $S_y(A) = 0$ for all $y \in B_1$, where $B_1 \subset B$ and $Q(B_1) = Q(B)$. Since $S.(A)$ is continuous, we get $S_y(A) = 0$ for all $y \in \overline{B_1}$. Since $\rho$ is strong, we have then $S_y(A) = 0$ for all $y \in \rho(\overline{B_1})$. But $\overline{B_1} \supseteq \rho(\overline{B_1}) \supset \rho(B_1) = \rho(B)$. □

CLAIM 3. *$P \otimes Q$ is not absolutely continuous with respect to $R$.*

PROOF. Let $H := \bigcup_{y \in Y}(C + y) \times \{y\}$. We claim that $H \in \mathfrak{A} \otimes \mathfrak{B}$. In fact, notice that $(x, y) \in H \Leftrightarrow x - y \in C$. Hence $H = g^{-1}(C)$, where $g(x, y) = x - y$. Since $S_y(H^y) = S_y(C + y) = 0$ for every $y \in Y$, we have $R(H) = 0$. On the other hand, $Q(H_x) = Q(-C + x) > 0$ for every $x \in X$, which gives $P \otimes Q(H) > 0$. □

A good example when the property (RF) does not take place is the following:

EXAMPLE 4.2. Let $\lambda$ be the Lebesgue measure on the real line. If $X = Y = [0, 1]$ are endowed with Lebesgue measurable sets and $\Delta$ is the diagonal of $[0, 1]^2$, then set $R(E) = \lambda(E \cap \Delta)/\sqrt{2}$ and $S_y = \delta_y$ (the measure concentrated in $\{y\}$) for all $y \in [0, 1]$. We have then $R([0, 1/2] \times [1/2, 1]) = 0$ but $S_{1/2}[0, 1/2] = 1$. According to Theorem 3.6 (or directly, by a simple calculation), there is a lifting $\pi \in \Lambda(\widehat{R})$ and there are liftings (uniquely determined in this case) $\sigma_y \in \Lambda(S_y)$ such that the equality (10) is satisfied:

$$[\pi(E)]^y = \sigma_y([\pi(E)]^y) \qquad \text{for all } y \in Y \text{ and } E \in \mathfrak{A} \widehat{\otimes}_R \mathfrak{B}.$$

As the measures $S_y$ are mutually singular, no rectangle formula holds true.

**5. An application to stochastic processes.** Let $\{\xi_y\}_{y \in Y}$ be an arbitrary real-valued stochastic process on $(X, \mathfrak{A}, P)$. If $\{\zeta_y\}_{y \in Y}$ is another stochastic process, then it is called $R$-*equivalent* to $\{\xi_y\}_{y \in Y}$ if, for each $y \in Y$, the equality $\xi_y = \zeta_y$ holds true a.e. $(S_y)$. $\{\zeta_y\}_{y \in Y}$ is then called a *modification* of $\{\xi_y\}_{y \in Y}$, and vice versa. $\{\xi_y\}_{y \in Y}$ is said to be $\widehat{R}$-*measurable* if the map $(x, y) \to \xi_y(x)$ is $\widehat{R}$-measurable. $\{\xi_y\}_{y \in Y}$ is bounded if $\sup_{y \in Y} \|\xi_y\|_{\mathfrak{L}^\infty(S_y)} < \infty$. There are several papers concerning the existence of $(P \widehat{\otimes} Q)$-measurable processes that are equivalent to a given process (cf. [3, 4, 16, 17]). If the initial process is already $(P \widehat{\otimes} Q)$-measurable, then one looks for its $(P \widehat{\otimes} Q)$-measurable modification behaving better than the original process. In general, a measurable process equivalent to a bounded $\{\xi_y\}_{y \in Y}$ is defined by setting $\zeta_y = \sigma(\xi_y)$, where $\sigma \in \Lambda(\widehat{P})$ and the initial process $\{\xi_y\}_{y \in Y}$ or the measure spaces satisfy some additional conditions. It is shown, however, in [3] that if the continuum hypothesis holds, then there exist nonpathological measure spaces $(X, \mathfrak{A}, P)$ and $(Y, \mathfrak{B}, Q)$, a lifting $\sigma \in \Lambda(\widehat{P})$ and an $(\mathfrak{A} \otimes \mathfrak{B})$-measurable stochastic process such that the lifting converts it into a non-$(P \widehat{\otimes} Q)$-measurable process (being a modification of the initial one). Thus, not always and not every lifting converts a measurable process into its measurable modification.

In the next theorem, we examine the problem of the existence of a measurable lifting modification of a measurable process in case of $R$ being not necessarily a product probability. We have not got any characterization of liftings converting measurable processes into their measurable modifications; we just show that if $\mathfrak{A}$ is not too large, then there exist liftings which always produce an $\widehat{R}$-measurable modification of an arbitrary $\widehat{R}$-measurable stochastic process. In fact, as the proof of Theorem 3.6 shows, in general there are a lot of such liftings. Since we assume only the existence of a product r.c.p., the measure $R$ may be quite far from the product measure $P \otimes Q$ and so the theorem describes a more general case than previously known results (when the separability of $\mathfrak{A}$ is assumed and no separability of the final process is needed).

THEOREM 5.1. *Assume that $\mathfrak{A}$ contains a countable algebra which is dense in $\mathfrak{A}$ with respect to $P$. Then for each bounded measurable stochastic process $\{\xi_y\}_{y \in Y}$ on $(X, \mathfrak{A}, P)$, there is a collection $\widetilde{\zeta} := \{\zeta_y\}_{y \in Y}$ of $\widehat{S}_y$-measurable functions $\zeta_y$ on $X$ and a collection of liftings $\sigma_y \in \Lambda(\widehat{S}_y)$, $y \in Y$, such that:*

(i) $\xi_y = \zeta_y$ *a.e.* $(S_y)$ *for all* $y \in Y$;
(ii) $\zeta_y = \sigma_y(\zeta_y)$ *for all* $y \in Y$;
(iii) *the map* $\widetilde{\zeta} : X \times Y \to (-\infty, +\infty)$ *is $\widehat{R}$-measurable.*

PROOF. In view of Theorem 3.6, there exist $\pi \in \Lambda(\widehat{R})$ and a family $\{\sigma_y \in \Lambda(\widehat{S}_y) : y \in Y\}$ such that, given process $\widetilde{\xi} = \{\xi_y\}_{y \in Y}$, we have

$$[\pi(\widetilde{\xi})]^y = \sigma_y([\pi(\widetilde{\xi})]^y) \quad \text{for all } y \in Y.$$

By Lemma 3.1, there exists $N_{\widetilde{\xi}} \in \mathfrak{B}_0$ such that

$$\xi_y = [\pi(\widetilde{\xi})]^y \quad \text{a.e. } (S_y) \text{ for all } y \notin N_{\widetilde{\xi}}.$$

We define now a collection $\widetilde{\zeta} := \{\zeta_y\}_{y \in Y}$ of $S_y$-measurable functions on $X$ by setting

$$\zeta_y = \sigma_y(\xi_y) \quad \text{for each } y \in Y.$$

Since $\pi(\widetilde{\xi})$ is $\widehat{R}$-measurable, one can easily see that $\{\zeta_y\}_{y \in Y}$ satisfies the required conditions. $\square$

**Acknowledgments.** The authors would like to thank the anonymous reviewers for their valuable remarks. Due to them, the results of the second section could have been essentially improved. Moreover, the presentation of the results is more readable now.

W. STRAUSS
FACHBEREICH MATHEMATIK
INSTITUT FÜR STOCHASTIK
  UND ANWENDUNGEN
ABTEILUNG FÜR FINANZ
  UND VERSICHERUNGSMATHEMATIK
UNIVERSITÄT STUTTGART
POSTFACH 80 11 40
D-70511 STUTTGART
GERMANY
E-MAIL: strauss@mathematik.uni-stuttgart.de

N. D. MACHERAS
DEPARTMENT OF STATISTICS
  AND INSURANCE SCIENCE
UNIVERSITY OF PIRAEUS
80 KAROLI AND DIMITRIOU STREET
185 34 PIRAEUS
GREECE
E-MAIL: macheras@unipi.gr

K. MUSIAŁ
INSTITUTE OF MATHEMATICS
WROCŁAW UNIVERSITY
PL. GRUNWALDZKI 2/4
50-384 WROCŁAW
POLAND
E-MAIL: musial@math.uni.wroc.pl